\documentclass[conf]{new-aiaa}
\usepackage[utf8]{inputenc}

\usepackage{graphicx}
\usepackage{amsmath}
\usepackage{multirow}
\usepackage[version=4]{mhchem}
\usepackage{siunitx}
\usepackage{svg}
\usepackage{longtable,tabularx}
\usepackage{colortbl}
\usepackage{geometry}
\usepackage{lscape}

\usepackage{pgf}
\usepackage{boldline}

\newcolumntype{x}[1]{>{\centering\arraybackslash\hspace{0pt}}p{#1}}
\def\cca#1{%
    \pgfmathsetmacro\calc{(#1-0.0)*100/(100.0)}%
    \edef\clrmacro{\noexpand\cellcolor{red!\calc}}%
    \clrmacro%
    \ifdim \calc pt>60pt\color{white}\fi{#1}%
}

\title{Industrial Application of Multidisciplinary Design Optimization with Uncertainties to a Pair of Telecommunication Satellites}

\author{O. Sapin \footnote{Research Engineer, Methods and Tools for the Development of Complex Systems, olivier.sapin@irt-saintexupery.com}}
\affil{IRT Saint Exupéry, Toulouse, France}

\author{L. Cousin \footnote{Research Engineer, Methods and Tools for the Development of Complex Systems, loic.cousin@irt-saintexupery.com}}
\affil{IRT Saint Exupéry, Toulouse, France}

\author{N. Roussouly \footnote{Research Engineer, Methods and Tools for the Development of Complex Systems, nicolas.roussouly@irt-saintexupery.com}}
\affil{IRT Saint Exupéry, Toulouse, France}

\author{A. Gazaix\footnote{Head of the MDO Competence Center, Methods and Tools for the Development of Complex Systems, anne.gazaix@irt-saintexupery.com}}
\affil{IRT Saint Exupéry, Toulouse, France}

\author{F. Gallard \footnote{Research Engineer, Methods and Tools for the Development of Complex Systems, françois.gallard@irt-saintexupery.com}}
\affil{IRT Saint Exupéry, Toulouse, France}

\author{M. De Lozzo \footnote{Research Engineer, Methods and Tools for the Development of Complex Systems, matthias.delozzo@irt-saintexupery.com}}
\affil{IRT Saint Exupéry, Toulouse, France}

\author{X. Fosse \footnote{End to end mission performance engineer, Space systems chief engineering, xavier.fosse@airbus.com}}
\affil{Airbus Defence and Space, Toulouse, France}

\author{N. Sarda \footnote{System Architecture Multi-Disciplinary Analysis and Optimization Expert, MBSE Design and Standards, nicolas.sarda@airbus.com}}
\affil{Airbus Defence and Space, Toulouse, France}

\author{G. Berthelin \footnote{Modeling and Simulation Architect, MBSE Design and Standards, gaspard.berthelin@airbus.com }}
\affil{Airbus Defence and Space, Toulouse, France}
\begin{document}

\maketitle

\begin{abstract}
In satellite design, it is common practice to add safety margins to the constraints to achieve conservative solutions that are robust to uncertainties.
This robustness often comes at the expense of performance and it may be more appropriate to include uncertainties in the definition of the design problem.
This work addresses such a challenge by applying techniques of multidisciplinary design optimization under uncertainty to an industrial use case.
The latter is a pair of telecommunication satellites launched together in a stacked configuration. Each satellite is a strongly coupled multidisciplinary system while the two satellites are not coupled at all.
This use case can be extended to an arbitrary number of satellites.
By considering uncertainty quantification techniques such as sensitivity analysis and reliability-based design optimization, 
this study demonstrates that accounting for uncertainties in the design problem results in a $66\%$ reduction in performance loss compared to the adding of a predefined safety margins,
while guaranteeing the feasibility of constraints with high probability.

\end{abstract}

\section{Introduction}\label{introduction}
\subsection{Context}\label{introduction:context}

Multidisciplinary Design Optimization (MDO) \cite{martins2013multidisciplinary} has yet to reach its full potential in industrial applications, despite its promise. The space domain, in particular, has been slow to adopt MDO, and more work is needed to make it more accessible and widely used. However, some progress has been made, with initiatives such as those undertaken at Airbus Defence and Space aiming to promote the use of MDO. \\

Uncertainties also pose a significant challenge in the context of MDO. Effective quantification and management (UQ\&M) of uncertainties are crucial, as uncertainties can have a major impact on the performance and feasibility of a design solution. Developing new methodologies and tools that can efficiently handle uncertainties will be essential to fully leverage the benefits of MDO and ensure its reliable application in various fields, including the space domain. \\

This paper aims to address these challenges by applying techniques of multidisciplinary design analysis and optimization (MDAO) under uncertainties (UMDAO) \cite{brevault2020aerospace} to the design of telecommunication satellites.
Section \ref{introduction:satellite_design} presents this design problem at a high level, 
while Section \ref{study_case} gives a detailed presentation of the MDO problems involved,
\textit{i.e.} objective, constraints, design variables, uncertainties and disciplines.
From a technological point of view,
Section \ref{introduction:rbdo} presents the notion of reliability-based MDO (RMDO) which is a suitable alternative to the use of safety margins.
Then, 
Section \ref{methodology} describes the methodology to set up this RMDO approach
and compare it with traditional approaches.
Section \ref{results} applies this methodology 
and shows that the RMDO approach leads to a more efficient solution than the use of safety margins.
Finally,
Section \ref{conclusion} gives conclusions and prospects for this work.

\subsection{Satellite design problem}\label{introduction:satellite_design}
\begin{figure}[hbt!]
\centering
\includegraphics[width=.9\textwidth]{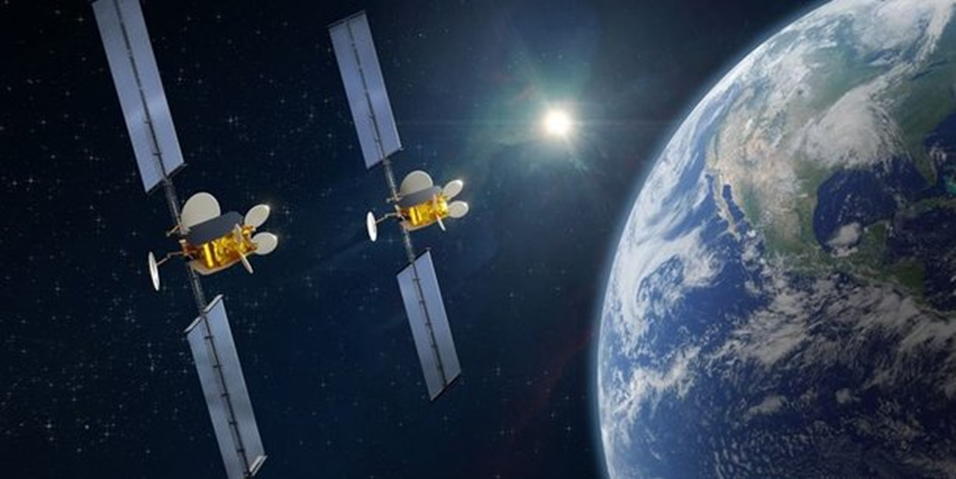}
\caption{Illustration of a dual launch of geostationary satellites sharing the same launch vehicle. © Copyright Airbus Defence and Space SAS 2021-2025}
\label{fig:pair_sat}
\end{figure}

The proposed engineering problem is the optimization of the orbital maneuver lifetime of a pair of geostationary satellites which are launched together, as in Figure~\ref{fig:pair_sat}.\\

The geostationary orbit (GEO) offers a unique vantage point allowing a satellite to be fixed relative to the Earth, providing a persistent view of the same ground target. It is typically used by communication satellites for broadcasting services. This orbit, however, is subject to small gravitational disturbances: the luni-solar attraction tends to move the satellite out of the equator plane, and the Earth non-spherical potential tends to move the satellite away from its operational longitude. These disturbances have to be compensated by regular station-keeping maneuvers, requiring a dedicated propulsion system and its associated propellant.\\

The orbital maneuver lifetime is the time a satellite can spend on its operational position given its available propellant for station-keeping. It is a very sensitive parameter, which directly drives the mission revenues, and its optimization is thus of utmost importance. 

Propellant is also required for the other phases: \\
\begin{itemize}
    \item the Electric Orbit Raising phase (EOR), starting from the launcher separation in the injection orbit to the geostationary orbit (GEO orbit). This phase spreads out from a few days to few months, depending on the strategy used to reach the GEO.\\
    \item the de-orbiting of the satellite at the end of its life. For GEO satellites, the strategy consists in sending the satellite into a graveyard orbit, at a specific higher altitude.\\
\end{itemize}

The computation of orbital maneuver lifetime shall also take into account diverse uncertainties (launcher performance, transfer and station-keeping delta-V, mass and CoM estimations, thruster performance...) so that the specified value for this critical parameter is reliable with a high confidence.

\subsection{Reliability‐based multidisciplinary design optimization}\label{introduction:rbdo}

Solving an MDO problem requires an MDO formulation, a.k.a. architecture \cite{martins2013multidisciplinary}. The latter consists in an optimization problem to be solved, including an objective function, constraints and design variables, together with a process organization that describes the execution sequence of disciplines and processes and the data exchange between these elements.
When uncertain variables affect disciplines, coupling variables, constraints and objective function become uncertain. In other words, for each value of these uncertain variables, the solution of the optimization problem differs and studying the impact of the uncertainty sources on the solution is a matter of interest. Furthermore, current industrial practices often manage uncertainties by generic safety margins, using a combination of worst‐case approaches, which may result in over‐conservative design solutions. Reliability-Based Design Optimization (RBDO) \cite{rbdo} then becomes an essential subject in order to provide more efficient and cost-effective solutions by explicitly accounting for uncertainty in the optimization problem, rather than relying on overly conservative safety margins.
Combining MDO formulation framework with RBDO leads to the context of RBMDO where dedicated methods have to be proposed to address the interdisciplinary couplings problem compared to classical RBDO. \\

In this paper,
this UMDO approach is based on the Multi-Disciplinary Feasible (MDF) formulation \cite{martins2013multidisciplinary}. Since the estimation of statistics by Monte Carlo (MC) sampling drastically increases the computational cost, 
it is rather proposed here to approximate the objective and constraints by Polynomial Chaos Expansion (PCE) surrogate models \cite{xiu2010numerical} and estimate the statistics of these surrogates by intensive MC sampling.

\section{The study case}\label{study_case}
\subsection{Disciplines}\label{study_case:description}
The study case focuses on the global propulsion system of two stacked GEO satellites. Several engineering disciplines are involved in the computation of the propellant budget and lifetime:\\
\begin{itemize}
    \item Launch vehicle performance: the injection orbit is dependent on the satellite mass, or in our case, of the stacked satellites. A higher launch mass means a lower energy orbit, typically with a lower apogee altitude, or a higher inclination relative to GEO, requiring a more costly transfer for the satellite, in terms of propellant and duration.\\
    \item Mission analysis: 
        \begin{itemize}
            \item Transfer strategy: the satellite performs maneuvers to reach GEO, consuming part of the available propellant. For an electric propulsion satellite, lower propellant consumption can be traded against higher transfer duration by applying a variable thrust strategy. This thrust strategy is defined by the relative electric orbit raising extension (noted below as $eor\_extension$). The higher this variable, the lower the consumption and the higher the duration of EOR.
            \item Station-keeping: the satellite performs regular maneuvers to maintain its operational position, but also occasional maneuvers for re-positioning and finally disposal on a graveyard orbit.\\
        \end{itemize}
    \item Mass \& centring: according to the propellant mass remaining in the tank after transfer, the Centre of Mass (CoM) variable position has an impact on the thrust direction of the robotic arms hosting the electric propulsion thrusters, which point to the CoM during station-keeping, thus changing the geometrical efficiency of the maneuvers, thus the station-keeping propellant consumption.\\
    \item Power budget: the satellite electric power is provided by its deployable solar arrays, which cells tend to degrade, notably during transfer when going through radiation belts. The transfer duration and trajectory has a significant impact on the power budget, which degrades the electric thruster performance.\\
    \item Electric thruster performance: the available thrust is dependent on the available power, but also on the propellant throughput experienced by individual thrusters due to aging effects.\\
\end{itemize}

The couplings between those disciplines makes this problem a classical MDAO use-case. The fact that 2 satellites or more are launched together, with potentially different characteristics such as dry mass and lifetime, adds another level of couplings.\\

\subsection{MDO without uncertainty}\label{study_case:deterministic}

The optimization problem, described in Equation \eqref{equation:deterministic}, aims to maximize the lifetime of the first satellite, which is both an objective and a design variable. Other design variables are also considered, such as:

\begin{itemize}
  \item the initial wet mass of the two satellites ($initial\_wet\_mass_0$ and $initial\_wet\_mass_1$), 
  \item the maximum accepted relative extension of electric orbit raising duration ($eor\_extension_0$ and $eor\_extension_1$), 
\end{itemize}
The design is constrained by: 
\begin{itemize}
  \item the launched mass, 
  \item the electric orbit raising duration for the two satellites, 
  \item the electric propellant consumption of the two satellites, 
  \item and the dry mass of each satellite.
\end{itemize}
The lifetime of the second satellite is considered fixed for the study. By convention, user data are noted $D_{some\_data}$, and come from specific configuration files.\\

\begin{equation}\label{equation:deterministic}
    \begin{aligned}
        \max_{
            \substack{
                initial\_wet\_mass_0\\
                initial\_wet\_mass_1\\
                eor\_extension_0\\
                eor\_extension_1\\
                lifetime_0
            }} \quad & lifetime_0 \\
        \textrm{s.t.} \quad & \forall i \in \{0, 1\}, \quad dry\_mass_i = D_{dry\_mass, i} \\
        & \forall i \in \{0, 1\}, \quad eor\_duration_i \leq D_{max\_eor\_duration, i} \\
        & \forall i \in \{0, 1\}, \quad elec\_propellant_i \leq D_{elec\_propellant\_max\_loading, i}\\
        & mass\_launched \leq D_{max\_mass\_launched}\\
    \end{aligned}
\end{equation}

\begin{figure}[hbt!]
\centering
\includegraphics[width=.9\textwidth]{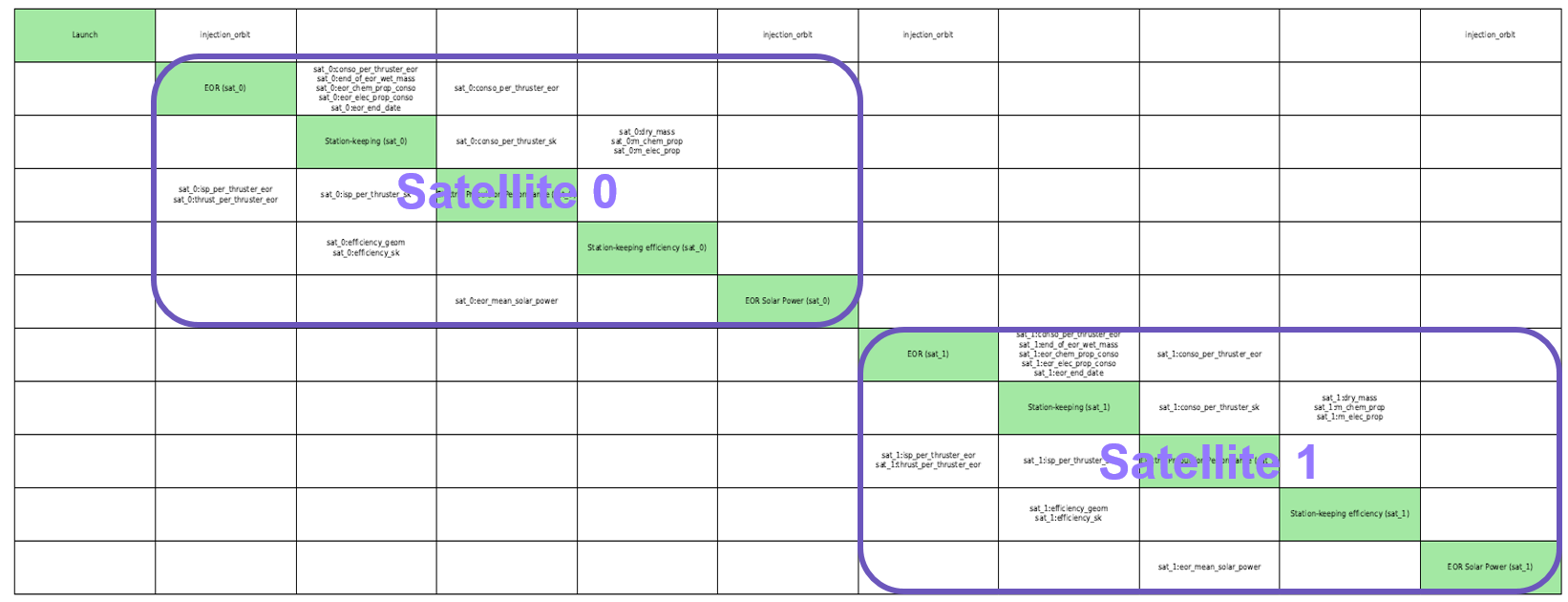}
\caption{N2 diagram representing all the disciplines and their couplings for a two-satellite configuration.}
\label{fig:N2}
\end{figure}

Figure \ref{fig:N2} represents the N2 diagram of the disciplines considered for the two-satellites configuration. 
We observe that the disciplines associated with a satellite include strong couplings, while the satellites are not coupled to each other but to the launch discipline.

Although this study on the two-satellites configuration, the implementation allows automatic extension to any number $n$ of satellites, as shown in the generalized Equation \eqref{equation:deterministic_n_sat}. 

\begin{equation}\label{equation:deterministic_n_sat}
    \begin{aligned}
        \max_{ 
            \substack{
                \begin{aligned}
                    {\scriptstyle \forall i \in \{0, ..., n-1\}} \quad & {\scriptstyle initial\_wet\_mass_i}\\
                    {\scriptstyle \forall i \in \{0, ..., n-1\}} \quad &  {\scriptstyle eor\_extension_i}\\
                    & {\scriptstyle lifetime_0}
                \end{aligned}
            }} \quad & lifetime_0 \\
        \textrm{s.t.} \quad & \forall i \in \{0, ..., n-1\}, \quad dry\_mass_i = D_{dry\_mass, i} \\
        & \forall i \in \{0, ..., n-1\}, \quad eor\_duration_i \leq D_{max\_eor\_duration, i} \\
        & \forall i \in \{0, ..., n-1\}, \quad elec\_propellant_i \leq D_{elec\_propellant\_max\_loading, i}\\
        & mass\_launched \leq D_{max\_mass\_launched}\\
    \end{aligned}
\end{equation}

\subsection{MDO with uncertainties}\label{study_case:uncertain}

The definition of the sources of uncertainty is one of the key steps to assess a reliable solution.
In this study, for each satellite, the sources of uncertainty in the system are modeled by nine normally distributed random variables (one design variable and eight parameters) with specified means and standard deviations:
\begin{itemize}
  \item the initial wet mass weighing,
  \item the electric propellant initial loading, 
  \item the electric propellant gauging accuracy at end of life,
  \item the specific impulse (ISP) dispersion in Electric Orbit Raising (EOR),
  \item the ISP dispersion in Station-Keeping, repositioning and re-orbiting,
  \item the efficiency dispersion in EOR,
  \item the efficiency dispersion in Station-Keeping, repositioning and re-orbiting,
  \item the delta-V (required total change in velocity) dispersion in EOR,
  \item the delta-V dispersion in Station-Keeping, repositioning and re-orbiting.
\end{itemize}
For the two-satellites configuration presented herein, a total of 18 uncertain variables were accordingly included in our model.\\

In the uncertain approach, the constraints of Equation \eqref{equation:deterministic} are considered to ensure that the probability of failure remains below a certain limit. This limit, denoted for simplification reasons in the following as $\varepsilon$, can be different for every constraint. This leads to the formulation of Equation \eqref{equation:rbdo}. \\

\begin{equation}\label{equation:rbdo}
    \begin{aligned}
        \max_{
            \substack{
                initial\_wet\_mass_0\\
                initial\_wet\_mass_1\\
                eor\_extension_0\\
                eor\_extension_1\\
                lifetime_0
            }} \quad & lifetime_0 \\
        \textrm{s.t.} \quad & \forall i \in \{0, 1\}, \quad \mathbb{E}[dry\_mass_i] = D_{dry\_mass, i} \\
        & \forall i \in \{0, 1\}, \quad \mathbb{P}[eor\_duration_i \geq D_{max\_eor\_duration, i}] \leq \varepsilon \\
        & \forall i \in \{0, 1\}, \quad \mathbb{P}[elec\_propellant_i \geq D_{elec\_propellant\_max\_loading, i}] \leq \varepsilon\\
        & \mathbb{P}[mass\_launched \geq D_{max\_mass\_launched}] \leq \varepsilon\\
    \end{aligned}
\end{equation}

Since the objective is treated as a design variable, it remains deterministic, whereas the constraints are subject to uncertainties. Consequently, the resulting optimization problem is a pure RBDO problem, consisting of minimizing a deterministic objective while ensuring that the probabilities of violating the constraints are small enough.

\section{Methodology}\label{methodology}

The whole study is implemented with the open source GEMSEO\footnote{GEMSEO (Generic Engine for Multidisciplinary Scenarios, Exploration and Optimization): \url{https://gemseo.org}.} \cite{gemseo_paper}\cite{UMDO2025} library, developed at IRT Saint~Exupéry.
The disciplines involved in the process are originally implemented in Scilab and are directly wrapped into GEMSEO through dedicated features. 

\subsection{MDO without uncertainty}\label{methodology;deterministic}

A first optimization and convergence analysis without uncertainty is performed in order to compute a deterministic reference point from which the reliability analysis is carried out by adding uncertainties.

We adopted the MDF formulation, where the coupling of strongly coupled disciplines is made consistent through a Multidisciplinary Design Analysis (MDA) at each iteration of the optimization process.

To perform the MDA, we utilized the \texttt{MDAChain} tool available in GEMSEO, which implements an advanced graph-based algorithm that examines the coupling graph and automatically identifies the strongly coupled disciplines.
This enables the decomposition of the full multidisciplinary system into weakly coupled multidisciplinary subsystems, composed of strongly coupled disciplines, that should be easier to solve.
In our specific case, the problem featured two strong couplings (one for each satellite), which we resolved using the nonlinear solver based on the fixed-point Jacobi method.\\

The optimization problem was solved using the COBYLA (Constrained Optimization BY Linear Approximations) algorithm from the NLopt library\footnote{NLopt: \url{https://nlopt.readthedocs.io/en/latest/NLopt_Algorithms}.}, which is well-suited for derivative-free optimization with nonlinear constraints, including both inequalities and equalities.

\subsection{Sensitivity analysis}\label{methodology:sensitivity_analysis}
Sensitivity analysis (SA) \cite{Iooss2015} is the study of how the model output variations can be explained by the different sources of uncertainty.
The design of complex systems often requires reliability assessments that involve a large number of uncertainties and low probability of failure estimations.
Estimating such rare event probabilities with crude MC sampling is computationally expensive, but the estimates are unbiased.
From this sampling, the SA can help to identify and filter out uncertain variables that have a low impact on variables of interest,
reducing the complexity of UMDO problems. 

In this study, the SA is performed by computing the total order Sobol' indices \cite{SOBOL2001271}. For each input and output of the model, this sensitivity index measures the full part of the variance explained by the input. They have been estimated with the Saltelli method \cite{SALTELLI2002280} implemented in OpenTURNS\footnote{(\url{https://openturns.github.io/openturns/latest/user_manual/_generated/openturns.SobolIndicesAlgorithm.html})} \cite{Baudin_2015}, which is wrapped in GEMSEO. To estimate these indices, a specific experimental design is needed, which involves generating two separate and independent samples. The number of evaluation required depends on the number of input variables. In our case, the estimation required $1000\times(18 + 2)=20 000$ evaluations of the multidisciplinary system, where 18 is the dimension of the uncertain space.

\subsection{Surrogate model using polynomial chaos expansions}
By leveraging surrogate model, also referred to as metamodel, UMDO problems can be tackled more efficiently as such a model can be evaluated for free. This approach allows the estimation of statistical moments, such as means and variances, as well as reliability measures, such as failure probabilities, to any desired precision by applying intensive MC sampling to the surrogate model.

A particularly effective surrogate modeling approach for uncertainty quantification and propagation is the PCE. The PCE of the output of interest $Y$ can be written as its truncated decomposition on an orthonormal polynomial basis:
\begin{equation*}
Y \approx c_0 + \sum_{i=1}^{k} c_i \phi_i(X)
\end{equation*}
with $X = (X_1, \dots, X_d)$, and where $\phi_i(X) = \prod_{j=1}^{d} \psi_{\tau_j(i),j}(X_j)$ is a polynomial of degree $p_i = \sum_{j=1}^{d} \tau_j(i)$ and  $\tau_j$ is an enumerating function mapping from $\mathbb{N}$ to $\mathbb{N}$.

The choice of the basis functions $\bigl(\psi_{i,j}\bigr)_{i,j \geq 0}$ depends on the input probability distributions. For Gaussian random variables, as it is the case here, Hermite polynomials are recommended. The coefficients $\bigl(c_i\bigr)_{i\geq0}$ can be determined by regression or quadrature. Here, we used a least-squares regression with a sparse strategy.

We used the OpenTURNS' PCE algorithm\footnote{(\url{https://openturns.github.io/openturns/latest/user_manual/_generated/openturns.FunctionalChaosAlgorithm.html})}, wrapped in GEMSEO, to generate the PCE, and  OpenTURNS' DOE algorithm optimal LHS\footnote{(\url{https://openturns.github.io/openturns/latest/user_manual/_generated/openturns.OptimalLHSExperiment.html})}, which is a Latin hypercube sampling (LHS) technique enhanced by simulated annealing. This global optimization technique starts from an initial LHS and improves it to maximize its discrepancy and so to get a better space-filling LHS.

\subsection{Reliability analysis}\label{methodology:reliability}

A reliability analysis was conducted to estimate the means, variances, and probabilities of constraint violation at the reference point, namely the optimum of the MDO problem. 
As explained above, PCEs were employed to alleviate the computational cost. They were constructed using a sample size of 40, which was selected to ensure its accuracy, and their quality was evaluated through cross-validation, resulting in an $\text{R}^2$ score greater than 0.99 for each constraint. The PCEs were built with a maximum degree of 2 and a limited number of terms (up to 30). Subsequently, the statistical estimates were obtained using 10 000 MC samples of the PCEs. \\

In a second step, generic margins were introduced on the thresholds of the active constraints to improve the reliability. A new  MDO was performed for each value of these margins, and a reliability analysis was conducted at these new reference points at the optimum.

\subsection{Reliablility based multidisciplinary design optimization}\label{methodology:rbdo}
\begin{figure}
    \centering
    \includegraphics[width=0.7\linewidth]{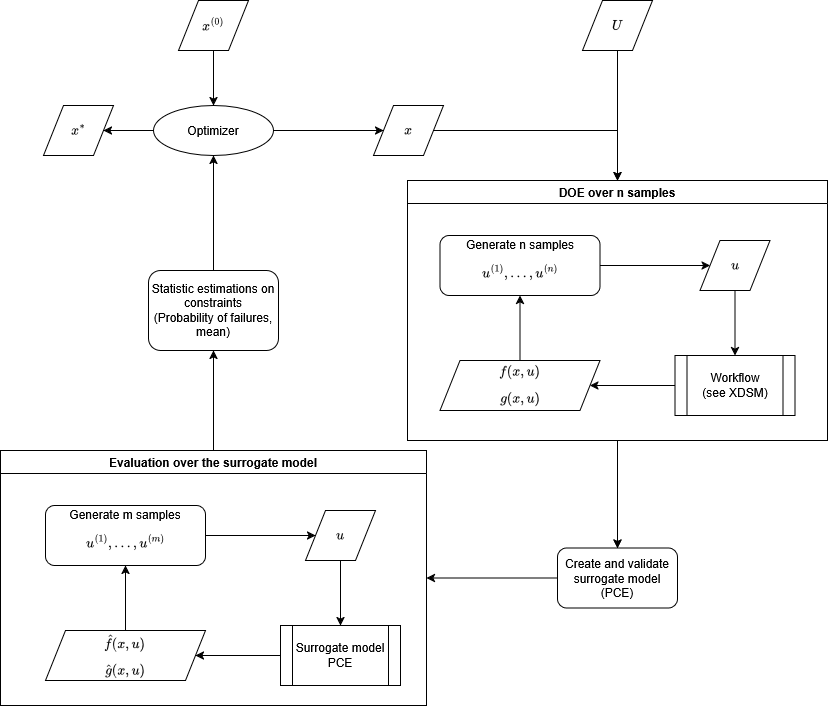}
    \caption{The RBMDO process.}
    \label{fig:methodology:rbdo}
\end{figure}

The UMDO problem was set up and solved using the GEMSEO-UMDO plugin\footnote{gemseo-umdo: \url{https://gemseo.gitlab.io/dev/gemseo-umdo}.}\cite{UMDO2025} using the strategy to generate PCE for all output quantities over the uncertain space at each iteration of the optimization loop and to use MC sampling to estimate the statistics, as shown in Figure \ref{fig:methodology:rbdo}.\\

We utilized the same setup to construct the PCEs at each optimization loop as the one that was employed for the reliability analysis, specifically a sample size of 40, a maximum degree of 2, and a limited number of terms (up to 30).
Their quality were monitored at each loop by logging the $\text{R}^2$ score, thereby detecting any potential poor approximation. 
Additionally, statistical estimates were obtained using 10 000 MC simulations at each iteration of the optimization loop.

\section{Results}\label{results}
\subsection{MDO without uncertainty}\label{results:deterministic}

Table \ref{tab:solution:deterministic} presents the value of the objective and constraints at the optimum of the MDO problem described by Equation \eqref{equation:deterministic}.
The units of constraints and their threshold values are deliberately not indicated in this document for confidentiality reasons.
We can see that this solution is feasible with two active constraints, namely $eor\_duration_0$ and $elec\_propellant_0$.  \\

\begin{table}[htp]
    \centering
    \renewcommand{\arraystretch}{1.6}
    \begin{tabular}{V{2}cV{2}cV{2}}
        \hlineB{2}
        \textbf{Output variable} & \textbf{Optimal value} \\
        \hlineB{2}
        $lifetime_0 (year)$ &  20.43 \\
        \hline
        $dry\_mass_0 - D_{dry\_mass, 0}$ & $- 3\times 10^{-7}$ \\
        \hline
        $dry\_mass_1 - D_{dry\_mass, 1}$ & $- 5\times 10^{-3}$ \\
        \hline
        $eor\_duration_0 - D_{max\_eor\_duration, 0}$ & $4\times 10^{-7}$ \\
        \hline
        $eor\_duration_1 - D_{max\_eor\_duration, 1}$ & $- 144$ \\
        \hline
        $elec\_propellant_0 - D_{elec\_propellant\_max\_loading, 0}$ & $2\times 10^{-7}$ \\
        \hline
        $elec\_propellant_1 - D_{elec\_propellant\_max\_loading, 1}$ & $- 100$ \\
        \hline
        $mass\_launched - D_{max\_mass\_launched}$ & $- 4246$ \\
        \hlineB{2}
    \end{tabular}
    \caption{Optimal solution of the MDO problem \eqref{equation:deterministic}.}
    \label{tab:solution:deterministic}
\end{table}

\begin{figure}
    \centering
    \includegraphics[width=0.8\linewidth]{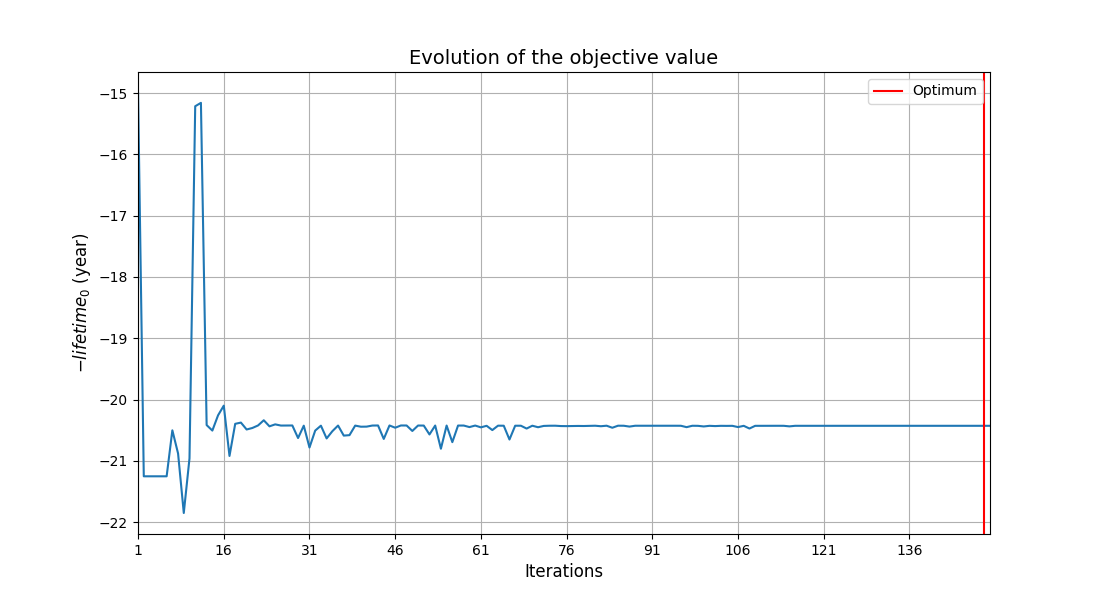}
    \caption{Evolution of the objective function.}
    \label{fig:deterministic:objective}
\end{figure}

Figure \ref{fig:deterministic:objective} shows the evolution of the objective during optimization. In our case, GEMSEO tries to minimize $- lifetime_0$.
Therefore, for each iteration shown on the x-axis, the corresponding negative objective value is plotted. The objective value starts to converge after 70 iterations.\\

\begin{figure}
    \centering
    \includegraphics[width=0.45\linewidth]{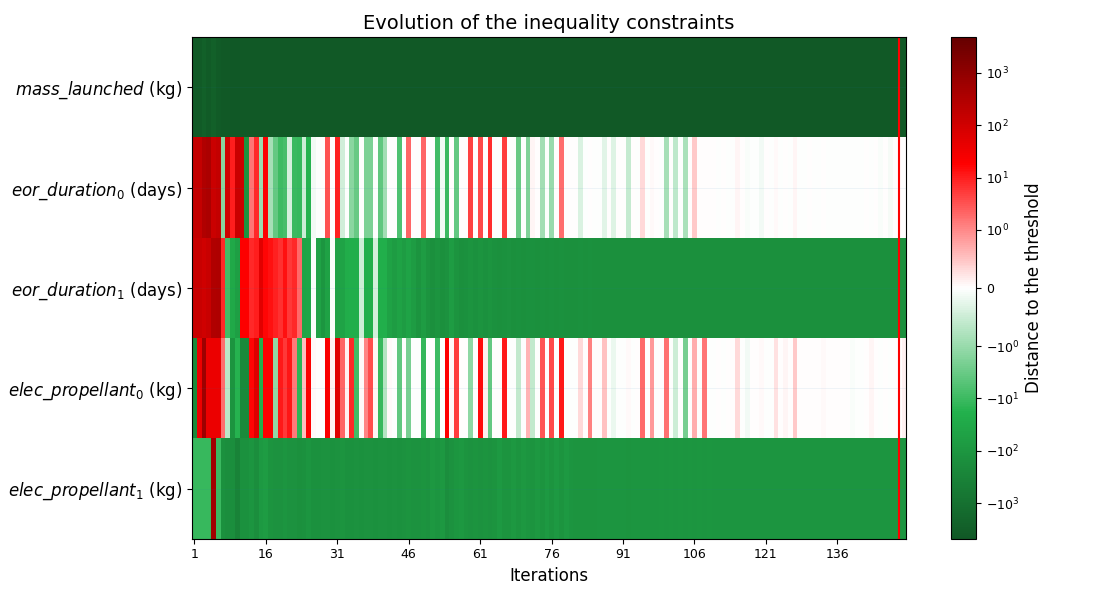}
    \includegraphics[width=0.45\linewidth]{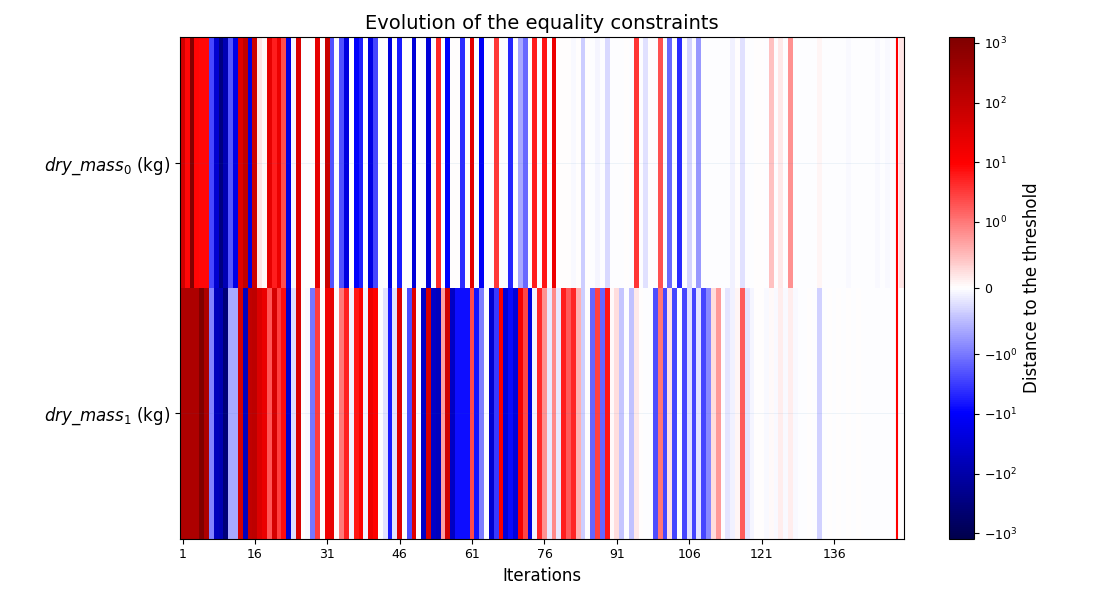}
    \caption{Evolution of inequality constraints (left) and equality constraints (right).}
    \label{fig:deterministic:constraints}
\end{figure}

Figure \ref{fig:deterministic:constraints} shows the evolution of the constraints during optimization. 
For each iteration shown on the x-axis, constraints are computed and the difference between their value and their threshold is plotted.
On the left graph, at a given iteration, if an inequality constraint is plotted in red, the system is not compliant with regards to it;
therefore, this constraint is said to be violated. If the color plotted is green, the system satisfies the requirements, and the constraint is inactive;
it is not constraining the system, so it has no influence on the system size. 
The white color appears once constraints are equal or close to zero, it means that the constraint is active:
if they are white at the optimum, they are sizing the system. Changing their threshold would change the system size.
The same explanation can be given for the right hand graph dealing with equality constraints. Violations of these constraints are indicated by red or blue colors, 
depending on whether the constraint value exceeds or falls below its threshold. In contrast, if an equality constraint is satisfied, it is represented in white.

\subsection{Sensitivity analysis}\label{results:uncertain:sensitivity_analysis}
A Sobol' analysis was performed from the deterministic solution. The total order Sobol' indices are given in Table \ref{table:results:sobol}. 
the redder the box, the greater the impact of the variable indicated in the column on the output indicated in the row.
The blue cells highlight the different uncertain parameters for which the cumulative total order Sobol' index do not overpass $10\%$. 
In other words, we point out the six variables that have a low impact on the constraints.

\begin{table}[htbp]
\centering
\renewcommand{\arraystretch}{1.6}
\setlength{\tabcolsep}{5pt}
\begin{tabular}{V{2}l*{18}{V{2}c}V{2}}
\hlineB{2}
 & \multicolumn{18}{cV{2}}{\textbf{Uncertain parameter}} \\
\clineB{2-19}{2}
\textbf{Constraint} & \rotatebox[origin=l]{90}{$u\_isp\_eor_0$} & \rotatebox[origin=l]{90}{$u\_isp\_eor_1$} &
\rotatebox[origin=l]{90}{$u\_eff\_eor_0$} & \rotatebox[origin=l]{90}{$u\_eff\_eor_1$} &
\cellcolor{blue!10}\rotatebox[origin=l]{90}{$u\_deltav\_eor_0$} &
\cellcolor{blue!10}\rotatebox[origin=l]{90}{$u\_deltav\_eor_1$} &
\rotatebox[origin=l]{90}{$u\_gauging_0$} & \rotatebox[origin=l]{90}{$u\_gauging_1$} &
\cellcolor{blue!10}\rotatebox[origin=l]{90}{$u\_loading_0$} &
\cellcolor{blue!10}\rotatebox[origin=l]{90}{$u\_loading_1$} &
\rotatebox[origin=l]{90}{$u\_isp\_sk_0$} & \rotatebox[origin=l]{90}{$u\_isp\_sk_1$} &
\rotatebox[origin=l]{90}{$u\_eff\_sk_0$} & \rotatebox[origin=l]{90}{$u\_eff\_sk_1$} &
\cellcolor{blue!10}\rotatebox[origin=l]{90}{$u\_deltav\_sk_0$} &
\cellcolor{blue!10}\rotatebox[origin=l]{90}{$u\_deltav\_sk_1$} &
\rotatebox[origin=l]{90}{$ u\_init\_wet\_mass_0 \quad$} & \rotatebox[origin=l]{90}{$u\_init\_wet\_mass_1 \quad$} \\
\hlineB{2}
$mass\_launched$        & \cca{0} & \cca{0} & \cca{0} & \cca{0} & \cca{0} & \cca{0} & \cca{0} & \cca{0} & \cca{0} &
\cca{0} & \cca{0} & \cca{0} & \cca{0} & \cca{0} & \cca{0} & \cca{0} & \cca{48}& \cca{54}\\
\hlineB{2}
$elec\_propellant_0$    & \cca{1} & \cca{0} & \cca{1} & \cca{0} & \cca{0} & \cca{0} & \cca{37}& \cca{0} & \cca{7} &
\cca{0} & \cca{22}& \cca{0} & \cca{23}& \cca{0} & \cca{3} & \cca{0} & \cca{0} & \cca{0} \\
\hlineB{2}
$dry\_mass_0$           & \cca{28}& \cca{0} & \cca{26}& \cca{0} & \cca{8} & \cca{0} & \cca{0} & \cca{0} & \cca{0} &
\cca{0} & \cca{16}& \cca{0} & \cca{17}& \cca{0} & \cca{1} & \cca{0} & \cca{6} & \cca{0} \\
\hlineB{2}
$eor\_duration_0$       & \cca{51}& \cca{0} & \cca{48}& \cca{0} & \cca{0} & \cca{0} & \cca{0} & \cca{0} & \cca{0} &
\cca{0} & \cca{0} & \cca{0} & \cca{0} & \cca{0} & \cca{0} & \cca{0} & \cca{0} & \cca{0} \\
\hlineB{2}
$elec\_propellant_1$    & \cca{0} & \cca{2} & \cca{0} & \cca{1} & \cca{0} & \cca{0} & \cca{0} & \cca{44}& \cca{0} &
\cca{9} & \cca{0} & \cca{20}& \cca{0} & \cca{19}& \cca{0} & \cca{2} & \cca{0} & \cca{0} \\
\hlineB{2}
$dry\_mass_1$           & \cca{0} & \cca{29}& \cca{0} & \cca{31}& \cca{0} & \cca{7} & \cca{0} & \cca{0} & \cca{0} &
\cca{0} & \cca{0} & \cca{13}& \cca{0} & \cca{13}& \cca{0} & \cca{2} & \cca{0} & \cca{4} \\
\hlineB{2}
$eor\_duration_1$       & \cca{0} & \cca{47}& \cca{0} & \cca{50}& \cca{0} & \cca{0} & \cca{0} & \cca{0} & \cca{0} &
\cca{0} & \cca{0} & \cca{0} & \cca{0} & \cca{0} & \cca{0} & \cca{0} & \cca{0} & \cca{0} \\
\hlineB{2}
\end{tabular}
\caption{Total order Sobol' indices in percent coming from 20 000 points.}
\label{table:results:sobol}
\end{table}

We investigated the impact of removing these variables on the solution of the UMDO problem. However, our results showed that eliminating these variables did not enable us to reduce the sample size required for constructing the PCE while maintaining the same level of accuracy, and it did not significantly decrease the computation time needed to build the PCE. Consequently, we retained all the $18$ variables in the RBMDO problem.

\subsection{Reliability analysis}\label{results:reliability}

The results of the initial reliability analysis are summarized Table \ref{table:results:failure_probabilities}, which presents the estimated failure probabilities and coefficients of variation for each constraint at the optimum solution of the MDO problem  \eqref{equation:deterministic}. As expected the active constraints display failure probabilities near $50\%$, which highlights the necessity of better accounting for the uncertainties affecting the constraints. \\

\begin{table}[htbp]
\centering
\renewcommand{\arraystretch}{1.6}
\setlength{\tabcolsep}{7pt}
\begin{tabular}{V{2}cV{2}cV{2}c|cV{2}}
\hlineB{2}
\multirow{2}{*}{\textbf{Constraint}} & \textbf{Coefficient of variation} & \multicolumn{2}{cV{2}}{\textbf{Reliability}} \\
\cline{3-4}
& $ = \frac{\text{standard deviation}}{\text{mean}}$ (\%) & Probability of failure (\%) & Confidence interval (\%) \\
\hlineB{2}
$dry\_mass_0$ & 0.22 & -- & -- \\
\hline
$dry\_mass_1$ & 0.24 & -- & -- \\
\hline
$eor\_duration_0$ & 1.28& 50.45 & 0.73 \\
\hline
$eor\_duration_1$ & 1.29 & 0.0 & 0.01 \\
\hline
$elec\_propellant_0$ & 1.43 & 50.30 & 0.73 \\
\hline
$elec\_propellant_1$ & 1.51 & 0.0 & 0.01 \\
\hline
$mass\_launched$ & 0.04 & 0.0 & 0.01 \\
\hlineB{2}
\end{tabular}
\caption{Failure probability and coefficient of variation for each constraint at the optimum solution of the MDO problem  \ref{equation:deterministic}.}
\label{table:results:failure_probabilities}
\end{table}

In a first rough approach, we investigate the addition of a margin on the threshold of the active constraint modifying the deterministic problem according to the Equation \eqref{equation:deterministic_with_margin},

\begin{equation}\label{equation:deterministic_with_margin}
    \begin{aligned}
        \max_{
            \substack{
                initial\_wet\_mass_0\\
                initial\_wet\_mass_1\\
                eor\_extension_0\\
                eor\_extension_1\\
                lifetime_0
            }} \quad & lifetime_0 \\
        \textrm{s.t.} \quad & \forall i \in \{0, 1\}, \quad dry\_mass_i = D_{dry\_mass, i} \\
        & eor\_duration_0 \leq D_{max\_eor\_duration, 0} - M_{eor\_duration_0}\\
        & eor\_duration_1 \leq D_{max\_eor\_duration, 1} \\
        & elec\_propellant_0 \leq D_{elec\_propellant\_max\_loading, 0} - M_{elec\_propellant_0}\\
        & elec\_propellant_1 \leq D_{elec\_propellant\_max\_loading, 1}\\
        & mass\_launched \leq D_{max\_mass\_launched}
    \end{aligned}
\end{equation}

\noindent We tested the following values for the margins using the standard deviations $ \sigma_{opt1} $ estimated at the optimum solution of the deterministic problem \ref{equation:deterministic}. 

\begin{equation*}
    \begin{aligned}
        M_{eor\_duration_0} & = k \times \sigma_{opt1}(eor\_duration_0)\\
        M_{elec\_propellant_0} & = k \times \sigma_{opt1}(elec\_propellant_0) \\
    \end{aligned}
\end{equation*}
with $ k = 1, 2, $ and $3$. The reliability analysis at the optimum solution of these deterministic problems presented and compared to the RBMDO solution Table \ref{table:results:comparison}.
We can see that the objective decreases slightly and that the probabilities of satisfying the constraints tend to 100\% when $k$ increases.

\subsection{Reliability Based Multidisciplinary Design Optimization}
The solution of the RBMDO problem given by Equation \eqref{equation:rbdo} is feasible, the value of the objective and the constraints statistics are presented Table \ref{equation:rbdo_res}.  The value of the objective function is reduced by $1.7\%$ compared to the solution of the MDO problem, which corresponds to a loss of approximately $4$ months in the lifetime. 
This was to be expected, as taking uncertainties into account in a design problem often leads to a more conservative average solution.\\

\begin{table}[htp]
    \centering
    \renewcommand{\arraystretch}{1.6}
    \begin{tabular}{V{2}cV{2}cV{2}}
        \hlineB{2}
        \textbf{RBMDO output} & \textbf{Optimal value} \\
        \hlineB{2}
        $lifetime_0$ (year) &  20.09 \\
        \hline
        $\mathbb{E}[dry\_mass_0] - D_{dry\_mass, 0}$ & $-7\times 10^{-5}$ \\
        \hline
        $\mathbb{E}[dry\_mass_1] - D_{dry\_mass, 1}$ & $1\times 10^{-5}$ \\
        \hline
        $\mathbb{P}[eor\_duration_0 \geq D_{max\_eor\_duration, 0}]$ & $0.$ \\
        \hline
        $\mathbb{P}[eor\_duration_1 \geq D_{max\_eor\_duration, 1}]$ & $0.$ \\
        \hline
        $\mathbb{P}[elec\_propellant_0 \geq D_{elec\_propellant\_max\_loading, 0}]$ & $3\times 10^{-3}$ \\
        \hline
        $\mathbb{P}[elec\_propellant_1 \geq D_{elec\_propellant\_max\_loading, 1}]$ & $0.$ \\
        \hline
        $\mathbb{P}[mass\_launched \geq D_{max\_mass\_launched}]$ & $0.$ \\
        \hlineB{2}
    \end{tabular}
    \caption{Optimal solution of the RBMDO problem (\ref{equation:rbdo}).}
    \label{equation:rbdo_res}
\end{table}

\noindent
Computational costs were significantly higher for the optimization under uncertainty, with each iteration taking roughly 100 times longer to complete than its deterministic counterpart, based on our computational infrastructure. This can be explained by the cost of obtaining the 40 training samples and building the sparse PCEs at each iteration of the optimization loop. \\

\begin{figure}
    \centering
    \includegraphics[width=0.8\linewidth]{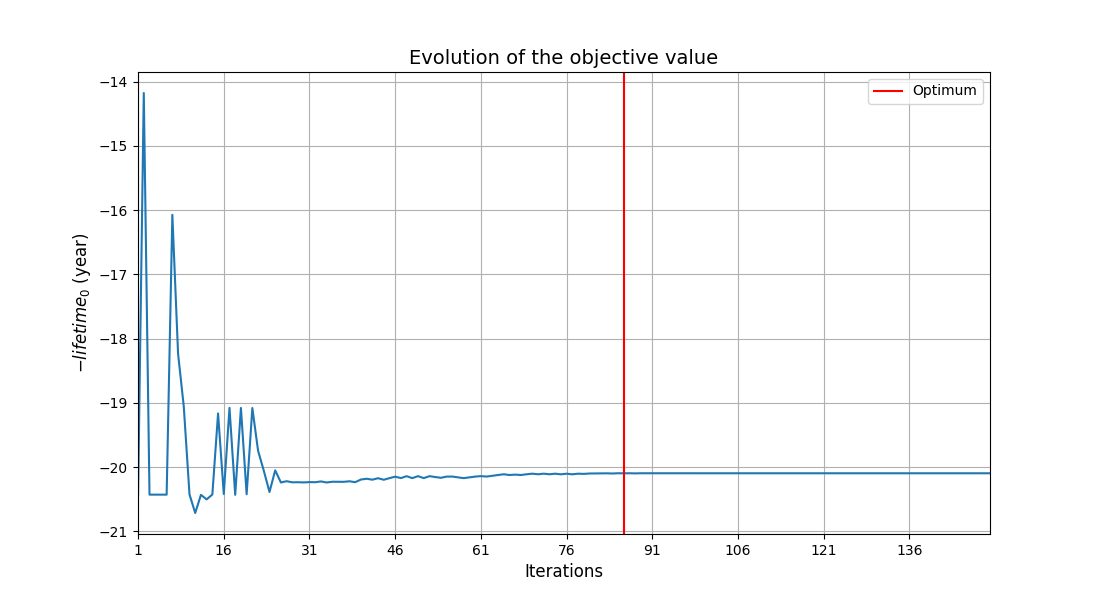}
    \caption{Evolution of the objective function during RBMDO.}
    \label{fig:rbdo:objective}
\end{figure}

\begin{figure}
    \centering
    \includegraphics[width=0.45\linewidth]{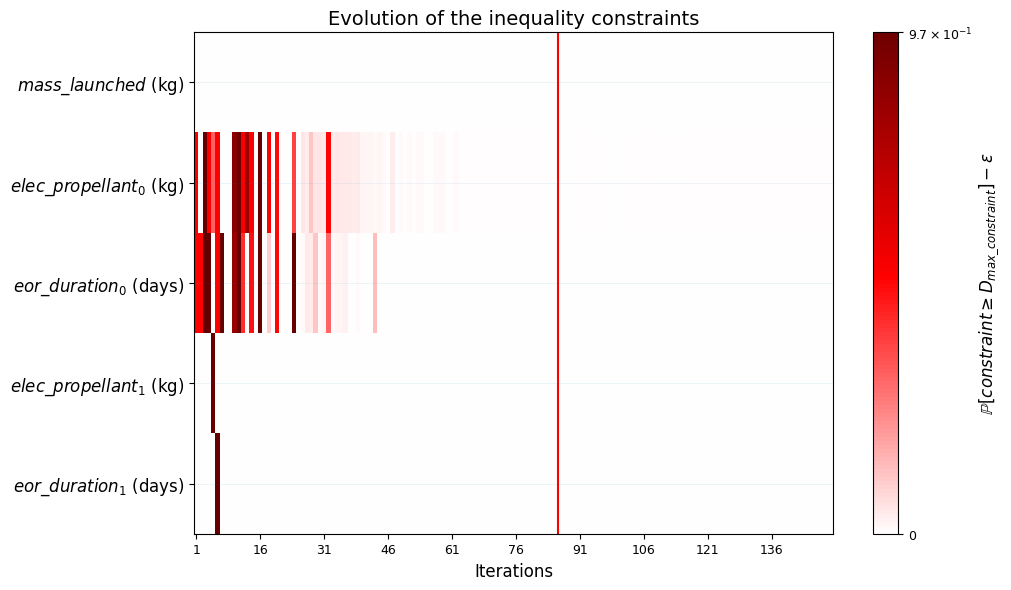}
    \includegraphics[width=0.45\linewidth]{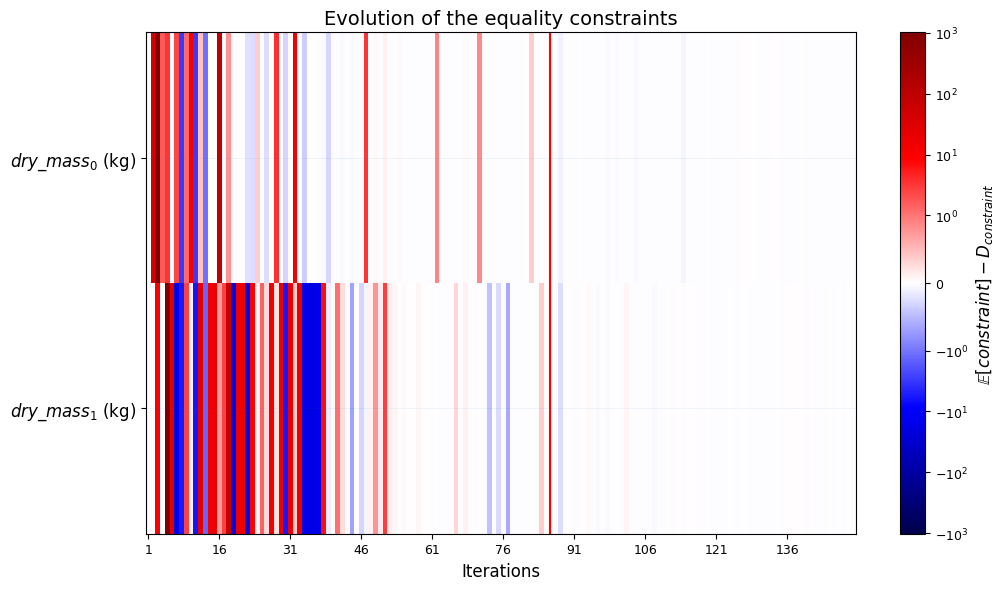}
    \caption{Evolution of inequality and equality constraints during RBMDO.}
    \label{fig:rbdo:constraints}
\end{figure}

\noindent
The results of the objective and constraints history are presented, respectively, in Figures \ref{fig:rbdo:objective} and \ref{fig:rbdo:constraints}.
For each iteration shown on the x-axis, the probability of failure of the inequality constraints and the expectation of the equality constraint are estimated
and the difference between their estimation and their threshold is represented.
On the left graph, the reliability threshold $\varepsilon$ is so small that we have not represented the negative value of the difference;
If the reliability condition is satisfied, it is represented in white, if not in red.
The same applies to the graph on the left, which deals with equality constraints. 
Violations of these constraints are indicated by red or blue colors, depending on whether the constraint average exceeds or falls below its threshold.
In contrast, if an equality constraint is satisfied, it is represented in white.

\begin{table}[htbp]
\centering
\renewcommand{\arraystretch}{1.6}
\setlength{\tabcolsep}{6pt}
\begin{tabular}{V{2}p{3.5cm}V{2}x{3cm}|x{3cm}V{2}x{3cm}V{2}}
\hlineB{2}
\multirow{2}{*}{\textbf{Method}} & \multicolumn{2}{cV{2}}{\textbf{Probability of failure (\%)}} & \textbf{Objective} \\
\cline{2-4}
& $eor\_duration_0$ & $elec\_propellant_0$ & $lifetime_0$ (year) \\
\hlineB{2}
MDO & 50.45 & 50.30 & 20.43 \\
\hline
RBMDO & 0.00 & 0.26 & 20.09 \\
\hline
MDO with $1\sigma$ margins & 17.70 & 16.39 & 20.08 \\
\hline
MDO with $2\sigma$ margins & 1.46 & 1.89 & 19.74 \\
\hline
MDO with $3\sigma$ margins & 0.00 & 0.01 & 19.41 \\
\hlineB{2}
\end{tabular}
\caption{Comparison of optimization methods with respect to failure probabilities and objective.}
\label{table:results:comparison}
\end{table}

Finally, we compare the objective values and failure probabilities for each method used in Table \ref{table:results:comparison}.
We observe that the strategy of introducing margins on the constraints leads to an unnecessary degradation of the objective value without providing an a priori guarantee of meeting the reliability conditions.
Specifically, the results indicate that a $3 \sigma$ margin would have been required to achieve the same level of reliability as in the RBMDO problem, resulting in a significant reduction of more than eight months in the lifetime compared to the RBMDO solution.\\

\section{Conclusion}\label{conclusion}

We presented a comprehensive application of multidisciplinary design optimization under uncertainty to the design of an space system involving a pair of telecommunication satellites launched in a stacked configuration.

By combining MDO and UQ\&M techniques using the GEMSEO open source software, including SA and RBDO, we have demonstrated the significant benefit of accounting for uncertainties in the design process. 

Our results showed that this approach can lead to more optimal and reliable design solutions compared to traditional deterministic methods, which rely on generic safety margins.
The study highlights the importance of incorporating uncertainty considerations into the design of complex systems, such as telecommunication satellites, and provides a valuable framework for future applications in the field.

Building on the multidisciplinary telecom system analyzed in this study, several optimization problems have been identified, providing a fertile ground for future investigation. A promising avenue for further research is to extend the current framework to accommodate constellations with more than two satellites, and to explore the optimization of orbital maneuver lifetime for all satellites using multi-objective optimization techniques, such as Pareto front analysis. This would enable the identification of optimal trade-offs between competing objectives, ultimately leading to improved design and operation of complex telecom systems.

\section*{Acknowledgment}
This work was supported by the French government under the France 2030 program under the reference ANR-10-AIRT-01.
The authors wish to acknowledge the R-Evol project members Airbus, Airbus Defence and Space, Liebherr, Capgemini, Expleo, Cenaero, CERFACS and INSA Toulouse for their support, financial funding and own knowledge.\\

The authors also thank J. Moulin for his invaluable contribution all along this work, his time and skills.

\bibliography{sample}

\begin{thebibliography}{10}
\newcommand{\enquote}[1]{``#1''}
\providecommand{\natexlab}[1]{#1}
\providecommand{\url}[1]{\texttt{#1}}
\providecommand{\urlprefix}{URL }
\expandafter\ifx\csname urlstyle\endcsname\relax
  \providecommand{\doi}[1]{\discretionary{}{}{}https://doi.org/#1}\else
  \providecommand{\doi}[1]{\discretionary{}{}{}\urlstyle{rm}\url{https://doi.org/#1}}\fi

\bibitem[{Martins and Lambe(2013)}]{martins2013multidisciplinary}
Martins, J.~R., and Lambe, A.~B., \enquote{Multidisciplinary design
  optimization: a survey of architectures,} \emph{AIAA journal}, Vol.~51,
  No.~9, 2013, pp. 2049--2075.

\bibitem[{Brevault et~al.(2020)Brevault, Balesdent, Morio
  et~al.}]{brevault2020aerospace}
Brevault, L., Balesdent, M., Morio, J., et~al., \emph{Aerospace system analysis
  and optimization in uncertainty}, Springer, 2020.

\bibitem[{Aoues and Chateauneuf(2010)}]{rbdo}
Aoues, Y., and Chateauneuf, A., \enquote{Benchmark study of numerical methods
  for reliability-based design optimization,} \emph{Structural and
  multidisciplinary optimization}, Vol.~41, No.~2, 2010, pp. 277--294.

\bibitem[{Xiu(2010)}]{xiu2010numerical}
Xiu, D., \emph{Numerical methods for stochastic computations: a spectral method
  approach}, Princeton university press, 2010.

\bibitem[{Gallard et~al.(2018)Gallard, Vanaret, Guénot, Gachelin, Lafage,
  Pauwels, Barjhoux, and Gazaix}]{gemseo_paper}
Gallard, F., Vanaret, C., Guénot, D., Gachelin, V., Lafage, R., Pauwels, B.,
  Barjhoux, P.-J., and Gazaix, A., \enquote{GEMS: A Python Library for
  Automation of Multidisciplinary Design Optimization Process Generation,}
  \emph{2018 AIAA/ASCE/AHS/ASC Structures, Structural Dynamics, and Materials
  Conference}, 2018.

\bibitem[{De~Lozzo et~al.(2025)De~Lozzo, Laboulfie, Sapin, Roussouly, Gallard,
  Aziz-Alaoui, Dechaume, and Gazaix}]{UMDO2025}
De~Lozzo, M., Laboulfie, C., Sapin, O., Roussouly, N., Gallard, F.,
  Aziz-Alaoui, A., Dechaume, A., and Gazaix, A., \enquote{Multidisciplinary
  design optimization under uncertainty: the open source capabilities of
  GEMSEO,} \emph{6th ECCOMAS Thematic Conference on Uncertainty Quantification
  in Computational Sciences and Engineering (UNCECOMP), Rhodes, Greece}, 2025.

\bibitem[{Iooss and Lema{\^i}tre(2015)}]{Iooss2015}
Iooss, B., and Lema{\^i}tre, P., \emph{A Review on Global Sensitivity Analysis
  Methods}, Springer US, Boston, MA, 2015, pp. 101--122.

\bibitem[{Sobol(2001)}]{SOBOL2001271}
Sobol, I., \enquote{Global sensitivity indices for nonlinear mathematical
  models and their Monte Carlo estimates,} \emph{Mathematics and Computers in
  Simulation}, Vol.~55, No.~1, 2001, pp. 271--280.
\newblock The Second IMACS Seminar on Monte Carlo Methods.

\bibitem[{Saltelli(2002)}]{SALTELLI2002280}
Saltelli, A., \enquote{Making best use of model evaluations to compute
  sensitivity indices,} \emph{Computer Physics Communications}, Vol. 145,
  No.~2, 2002, pp. 280--297.

\bibitem[{Baudin et~al.(2015)Baudin, Dutfoy, Iooss, and Popelin}]{Baudin_2015}
Baudin, M., Dutfoy, A., Iooss, B., and Popelin, A.-L., \emph{OpenTURNS: An
  Industrial Software for Uncertainty Quantification in Simulation}, Springer
  International Publishing, 2015, p. 1–38.

\end{thebibliography}

\end{document}